\documentclass[11pt, oneside, a4paper]{article}
\usepackage{graphicx}
\usepackage{tikz}
\usepackage{cite}
\usepackage{hyperref}
\usepackage{amsfonts}
\usepackage{fancyhdr,amsthm,amsmath,amssymb,mathrsfs}
\usepackage{bookmark}
\usetikzlibrary{intersections}

\newtheorem{thm}{\bf Theorem}[section]
\newtheorem{lem}[thm]{Lemma}
\newtheorem{prop}[thm]{Proposition}
\newtheorem{cor}[thm]{Corollary}

\theoremstyle{definition}

\newtheorem{eg}{Example}
\numberwithin{equation}{section}
\newenvironment{prf}{\noindent {\bf Proof.}\rm}{\qed}
\newcommand{\pth}{{\rm Path}}

\newcommand{\cnw}{{\rm C}(W)}
\newcommand{\cnh}[1]{{\rm C}(#1)}

\newcommand{\In}{{\rm Inv}(\Gamma)}

\newcommand{\al}{\alpha}

\newcommand{\G}{\Gamma}

\begin{document}
\title{\Large\bf On lattice of congruences on graph inverse semigroups\thanks{Partially supported by Chongqing Natural Science
Foundation (cstc2019jcyj-msxmX0435) and the Fundamental Research Funds for the Central Universities (XDJK2016B038).}}
\author{Yongle Luo and Zhengpan Wang\thanks{Corresponding author} \\
{\small\em School of Mathematics and Statistics, Southwest University, $400715$ Chongqing, China} \\
{\small\em E-mail: lechoke@qq.com, \; zpwang@swu.edu.cn} }
\date{}
\maketitle

\begin{abstract}
Congruences on a graph inverse semigroup were recently described in terms of the underline graph. Based on such descriptions, we show that the lattice of congruences on a graph inverse semigroup is upper semimodular but not lower semimodular.
\medskip


\end{abstract}

\section{Introduction and main result}

Graph inverse semigroups were first introduced by Ash and Hall in 1975 for graphs which have at most one edge from one vertex to another (see \cite{AshHall}). They gave necessary and sufficient conditions for such graph inverse semigroups to be congruence free. The general definition of graph inverse semigroups can be found in \cite{Jones, JonesLawson, MesyanMitchell} etc. These inverse semigroups are related to various types of algebras, for instance, $C^\ast$-algebras, Cohn path algebras, Leavitt path algebras and so on (see \cite{AbramsPino, AlaliGilbert, JonesLawson, Krieger, KumjianPaskRaeburnRenault, KumjianPaskRaeburn, Paterson}). Graph inverse semigroups have also been studied in their own right in recent years (for instance, see \cite{AlaliGilbert, ChihPlessas, Jones, JonesLawson, MesyanMitchell, MesyanMitchellMoraynePeresse, Wang}). Congruences on graph inverse semigroups were particularly discussed in \cite{MesyanMitchell, Wang}. The second author Wang built a one-to-one correspondence of the set of congruences on graph inverse semigroups onto the set of so-called congruence triples of graphs (see \cite{Wang}). Based on such a description, we prove in this paper that the lattice of congruences on a graph inverse semigroup is upper semimodular but not lower semimodular.

A {\it directed graph} $\G=(V, E, r, s)$ consists of $V, E$ and functions $r, s : E \rightarrow V$. The elements of $V$ are called {\it vertices} and the elements of $E$ are called {\it edges}. For each edge $e$, $s(e)$ is the {\it source} of $e$ and $r(e)$ is the {\it range} of $e$. For any set $X$, we denote the {\it cardinality} of $X$ as $|X|$, and for any vertex $v$, $|s^{-1}(v)|$ is called the {\it index} of $v$. In this paper, we always denote a directed graph $\G=(V, E, r, s)$ simply as $\G$.

The {\it graph inverse semigroup} $\In$ of $\G$ is the semigroup with zero generated by $V$ and $E$, together with the set $E^\ast$, of variables $\{e^\ast : e \in E\}$, satisfying the following relations for all $u, v \in V$ and $e, f \in E$:
\begin{description}
\item[{\rm(GI1)}] $u v = \delta_{u,v} u$;
\item[{\rm(GI2)}] $s(e) e = e r(e) = e$;
\item[{\rm(GI3)}] $r(e) e^\ast = e^\ast s(e) = e^\ast$;
\item[{\rm(CK1)}] $f^\ast e = \delta_{f,e} r(e)$.
\end{description}
A {\it path} in $\G$ is a sequence $\al := e_1 \cdots e_n$ of edges such that $r(e_i) = s(e_{i+1})$ for $i = 1, \cdots, n-1$. In such a case, $s(\al) := s(e_1)$ is the {\it source} of $\al$, $r(\al) := r(e_n)$ is the {\it range} of $\al$.  The path $\al$ is called a {\it cycle} if $s(\al) = r(\al)$ and $s(e_i) \neq s(e_j)$ whenever $i \neq j$. Denote the set of all paths in $\G$ as $\pth(\G)$, and for an arbitrary path $\alpha= e_1 e_2...e_n$, denote the set $\{s(e_i): i=1, 2,...,n\}$ as $u(\alpha)$. In particular, we also view a vertex $v$ of $V$ as a path, $s(v) = r(v) = v$, and $u(v)= \emptyset$. For any $V_1 \subseteq V$, ${\cnh{V_1}}$ is the set of all cycles $c$ with $u(c) \subseteq V_1$.

On graph $\G$, a subset $H$ of $V$ is {\it hereditary} if $s(e) \in H$ always implies $r(e) \in H$ for all $e \in E$. Let $H$ be a hereditary subset of $V$. For any $v \in V$, we call $|\{e \in E:  s(e) = v, r(e) \notin H \}|$ the {\it index of $v$ relative to} $H$. Let $W$ be a subset of the set of vertices with index 1 relative to $H$. Clearly, $H$ and $W$ are two disjoint subsets. We call $f: \cnw \rightarrow \mathbb{Z}^+ \cup \{\infty\}$ a {\it cycle function on $W$}, if $f$ is invariant under cyclic permutations. A {\it congruence triple of $\G$}, written by $(H, W, f)$, consists of a hereditary subset $H$, a subset $W$ of the set of vertices with index 1 relative to $H$ and a cycle function $f$ on $W$ (see \cite{Wang}). Given a congruence triple $(H, W, f)$, for convenience, we extend the cycle function $f$ to $\bar{f} : \cnh{V} \rightarrow \mathbb{Z}^+ \cup \{\infty\}$ with respect to $H$ and $W$ as what follows: for any $c \in \cnh{V}$,
\begin{equation} \label{extending}
\bar{f}(c) = {\begin{cases}
1 & \mbox{if } c \in \cnh{H}, \\
f(c)  & \mbox{if } c \in \cnw, \\
\infty & \mbox{otherwise}.
\end{cases}}
\end{equation}
Obviously, these extended cycle functions one-to-one correspond to original cycle functions. In what follows, we will use extended cycle functions to replace cycle functions in a congruence triple and still use $(H, W, f)$ instead of $(H, W, \bar{f})$. Particularly, if $\cnh{V} = \emptyset$, then we denote a corresponding congruence triple by $(H, W, \emptyset)$. We denote the set of all replaced congruence triples of $\G$ as $\mathcal{CT}(\G)$.

For any $m_1, m_2 \in \mathbb{Z}^+ \cup \{\infty\}$, we denote $(m_1, m_2)$ their {\it greatest common divisor} and $[m_1, m_2]$ their {\it least common multiple}. We naturally assume that any element in $\mathbb{Z}^+ \cup \{\infty\}$ divides $\infty$, which means that, for all $m\in \mathbb{Z}^+ \cup \{\infty\}$, we have $(m, \infty)= m$ and $[m, \infty]= \infty$. Now we define a binary relation on $\mathcal{CT}(\G)$: for any $(H_1, W_1, f_1), (H_2, W_2, f_2) \in \mathcal{CT}(\G)$, $(H_1, W_1, f_1) \leq (H_2, W_2, f_2)$ if $H_1 \subseteq H_2$, $W_1 \setminus H_2 \subseteq W_2$ and $f_2(c) \mid f_1(c)$ for any $c \in \cnh{V}$. The next corollary illustrates that this binary relation coincides with the partial order on the set of all congruence triples defined in \cite{Wang}.

\begin{cor}
Given $(H_1, W_1, f_1), (H_2, W_2, f_2) \in \mathcal{CT}(\G)$ with $H_1 \subseteq H_2$ and $W_1 \setminus H_2 \subseteq W_2$, then $f_2(c)\mid f_1(c)$ for any $c \in \cnh{V}$ if and only if $f_2(c)\mid f_1(c)$ for any $c \in{\cnh{W_1}\cap \cnh{W_2}}$.
\end{cor}

\begin{prf}
The necessity part is clear. Note that ${\cnh{W_1}\cap \cnh{W_2}} = \cnh {W_1 \cap {W_2}}$. For any $c\in \cnh{V}$, if $c \in {\cnh{H_2}}$, then we certainly obtain that $f_2(c) \mid f_1(c)$ since $f_2(c)=1$ by (\ref{extending}); if $c\in {\cnh{V} \setminus (\cnh{W_1 \cap W_2)}\cup \cnh{H_2})}$, then we also obtain that $f_2(c)\mid f_1(c)$ since $f_1(c)=\infty$.
\end{prf}
\medskip

Let $\rho$ be a congruence on $\In$. We define $H = 0 \rho \cap V$ and $W = \{s(e) \in V \setminus H: ee^\ast \rho s(e) \mbox{ for } r(e) \notin H\}$. For $c \in {\cnh{W}}$, let $f(c)$ be the smallest positive integer $m$ such that $c^m \rho s(c)$. If no power of $c$ is equivalent to $s(c)$, then define $f(c) = \infty$. Thus, $T(\rho) = (H, W, f)$ is a congruence triple of $\G$. Conversely, let $(H, W, f)$ be a congruence triple of $\Gamma$ and let $\wp(H, W, f)$ denote the congruence generated by all pairs $(v, 0)$ for $v \in H$, $(ee^\ast, w)$ for $w \in W$ with $s(e) = w$ and $r(e) \notin H$, and $(c^{f(c)}, s(c))$ for $c \in \cnw$ with $f(c) \in \mathbb{Z}^+$. Then $\wp(H, W, f)$ is a congruence on $\In$. Denote the set of all congruences on the inverse semigroup $\In$ as $\mathcal{C(S)}$. Wang proved in \cite{Wang} that functions $\wp$ and $T$ between $\mathcal{C}(S)$ and $\mathcal{CT}(\G)$ are inverses of each other. In fact, we also have

\begin{prop} \label{isom}
The partial ordered sets $(\mathcal{CT}(\G), \leq )$ and $(\mathcal{C}(S), \subseteq)$ are order isomorphic.
\end{prop}

\begin{prf}
For all $\rho_1, \rho_2 \in \mathcal{C}(S)$ with $\rho_1 \subseteq \rho_2$, if $T(\rho_1) = (H_1, W_1, f_1), T(\rho_2) = (H_2, W_2, f_2)$, then it is proved at the end of Section~1 in \cite{Wang} that $(H_1, W_1, f_1)\leq (H_2, W_2, f_2)$. That is to say, $T$ is order preserving.

We now prove that $\wp$ is order preserving as well. Given $(H_1, W_1, f_1), (H_2, W_2, f_2) \in \mathcal{CT}(\G)$ with $(H_1, W_1, f_1)\leq (H_2, W_2, f_2)$, we write congruences $\rho_1 = \wp(H_1, W_1, f_1)$, $\rho_2 = \wp(H_2, W_2, f_2)$ respectively generated by relations $R_1, R_2$, where
\begin{align*}
R_{i}= & \{(v, 0): v\in{H_i}\} \cup  \{(ee^\ast, w): w\in{W_i}, e\in{s^{-1}(w)}, r(e)\notin{H_i}\}\\
&\cup \{(c^{f_i(c)}, s(c)): c\in{\cnh{W_i}}, f_i (c)\in {\mathbb{Z}^{+}} \}, ~~~i= 1, 2.
\end{align*}
To show $\rho_1 \subseteq \rho_2$, we only need to prove that $R_1 \subseteq \rho_2$. For every $v \in{H_1}$ we observe from $ H_1 \subseteq H_2$ that $(v,0)\in{R_2}$. For any $w\in{W_1}$ and $e\in{s^{-1}(w)}$ with $r(e)\notin{H_1}$, if $ w\notin{H_2}$, then we obtain that $w\in{W_1\setminus{H_2}} \subseteq{W_2}$ which means that $e$ must be the unique edge with $s(e) = w$ and $r(e) \notin{H_2}$. That is, $(ee^\ast,w) \in R_2$. If $w \in{H_2}$, then it follows from the definition of a hereditary subset that $r(e) \in H_2$ which leads to $(r(e), 0) \in R_2$. So we get $(e e^\ast, 0) \in \rho_2$. Note that $(w, 0) \in R_2$, then we have $(e e^\ast, w) \in \rho_2$. Take $c \in \cnh{W_1}$ with $f_1(c) \in \mathbb{Z}^{+}$. If $s(c)\notin{H_2}$, then $s(c)\in{W_1\setminus{H_2}} \subseteq W_2 $. It is clear that a vertex from $W_1$ lies either in $W_2$ or $H_2$, so all vertices in $u(c)$ lie in $W_2$. Thus, $f_2 (c)\mid f_1 (c)$ implies that $f_2(c) \in \mathbb{Z}^+$. It follows from $(c^ {f_2 (c)}, s(c))\in R_2$ that $(c^{f_1 (c)}, s(c))\in \rho_2$. If $s(c)\in H_2$, then $(s(c), 0)\in R_2$ which means that $(c^{f_1 (c)}, 0)\in \rho_2$. Hence $(c^ {f_1 (c)}, s(c))\in \rho_2$. In conclusion, we get $R_1 \subseteq \rho_2$.
\end{prf} \medskip

Consequently, $(\mathcal{CT}(\G), \leq )$ is a lattice. On a lattice $L$, we denote the {\it meet} and the {\it join} of $a$ and $b$ respectively as $a \wedge b $ and $a \vee b $. We say that $a$ {\it covers} $b$, and write $a\succ b$ or $b \prec a$, if $a \gneq b$ and if there is no $x$ in $L$ such that $a \gneq x \gneq b$. $L$ is said to be {\it distributive} if, for all $a, b, c $ in $L$, $ (a\vee b)\wedge c =(a\wedge c)\vee (b\wedge c)$. $L$ is said to be {\it modular} if, for all $a, b, c $ in $L$, $a\leq c$ implies that $(a\vee b)\wedge c = a\vee (b\wedge c).$  $L$ is said to be {\it upper semimodular} if, for all $a, b$ in $L$, $a\succ a\wedge b, b\succ a\wedge b$ imply $a\vee b \succ a, a\vee b \succ b$. Dually, $L$ is said to be {\it lower semimodular} if, for all $a, b$ in $L$, $a\vee b \succ a, a\vee b \succ b$ imply $a\succ a\wedge b, b\succ a\wedge b$. It is known that a modular lattice is both upper and lower semimodular (see \cite{Howie}).

We know that the lattice of normal subgroups of a group is modular. However, the lattice of congruences on an inverse semigroup is not even semimodular. For graph inverse semigroups, we have

\begin{thm} \label{2 upper semimodular}
The lattice $(\mathcal{CT}(\G), \leq)$ is upper semimodular, but not necessarily lower semimodular.
\end{thm}

\section{Proof of Theorem~\ref{2 upper semimodular}}

First we give some lemmas on covers of congruence triples.

\begin{lem} \label{cover of H}
For any $(H_0, W_0, f_0), (H, W, f) \in \mathcal{CT}(\G)$ with $(H_0, W_0, f_0) \prec (H, W, f)$ and $H_0 \subsetneq H$, the following statements hold.
\begin{description}
\item[{\rm(1)}] $W_0\setminus H =W$ and $f_0(c)= f(c)$ for any $c \in \cnh{W_0}$.
\item[{\rm(2)}] There exists no vertex in $H \setminus {(H_0 \cup W_0})$ which has index $1$ relative to $H_0$.
\end{description}
\end{lem}

\begin{prf}
$(1)$ Note from $(H_0, W_0, f_0) \prec (H, W, f)$ that $W_0 \setminus H \subseteq W$. Define a cycle function $g$ with respect to $H$ and $W_0 \setminus H$ by $g(c) = f_0(c)$ for any $c \in \cnh{W_0 \setminus H}$. Then we see that $(H, W_0 \setminus H, g) \in \mathcal{CT}(\G)$ and also that $(H_0, W_0, f_0) \lneq (H, W_0 \setminus H, g) \leq (H, W, f)$. So we have $W_0 \setminus H = W$.

Defining a cycle function $f'$ with respect to $H_0$ and $W_0$ by $f'(c)= f(c)$ for any $c\in \cnh{W_0}$. Then we have $(H_0, W_0, f_0) \leq (H_0, W_0, f') \lneq (H, W, f)$ which leads to $f_0= f'$. Thus, we obtain that $f_0(c) =f(c)$ for any $c\in \cnh{W_0}$.

$(2)$ Suppose that there exists a vertex $v \in {H\setminus {(H_0 \cup W_0)}}$ such that $v$ has index $1$ relative to $H_0$. Then it is clear that $(H_0, W_0 \cup {\{v\}}, f_0) \in \mathcal{CT}(\G)$ and that $(H_0, W_0, f_0) \lneq (H_0, W_0 \cup {\{v\}}, f_0) \lneq (H, W, f)$ which contradicts the hypothesis.
\end{prf}
\medskip

For any $(H_0, W_0, f_0), (H, W, f) \in \mathcal{CT}(\G)$, if $(H_0, W_0, f_0) \prec (H, W, f)$, then we remark either from Lemma~\ref{cover of H} that $W_0 \setminus H = W$ which means that $W \subseteq W_0$ when $H_0 \subsetneq H$ or from the definition of the partial order on $\mathcal{CT}(\G)$ that $W_0 \subseteq W$ when $H_0 = H$. Then the following corollary is immediate.

\begin{cor} \label{w0subsetneqw}
For any $(H_0, W_0, f_0), (H, W, f) \in \mathcal{CT}(\G)$ with $(H_0, W_0, f_0) \prec (H, W, f)$, if $W_0 \subsetneq W$, then $H_0 = H$.
\end{cor}

\begin{lem} \label{cover of W}
For any $(H_0, W_0, f_0), (H, W, f) \in \mathcal{CT}(\G)$ with $W_0 \subsetneq W$, $(H_0, W_0, f_0)\prec(H, W, f)$ if and only if $H_0 = H$, $f_0 = f$ and $|W\setminus W_0|=1$.
\end{lem}

\begin{prf}
If $(H_0, W_0, f_0)\prec(H, W, f)$, then we see from Corollary~\ref{w0subsetneqw} that $H_0 = H$ since $W_0 \subsetneq W$. Thus, we obtain that $(H, W, f_0) \in \mathcal{CT}(\G)$ and that
$(H_0, W_0, f_0) \lneq (H, W, f_0) \leq (H, W, f)$ so
that $f_0 = f$. Moreover, for any $W'$ such that $W_0 \subsetneq W' \subseteq W$, we observe that $(H, W', f_0) \in \mathcal{CT}(\G)$. So by $(H_0, W_0, f_0) \lneq (H, W', f_0) \leq (H, W, f)$, we obtain that $|W\setminus W_0|=1$. The converse part is obvious.
\end{prf}

\begin{lem} \label{samehw}
For any $(H, W, f_0), (H, W, f)\in \mathcal{CT}(\G)$, $(H, W, f_0) \prec (H, W, f)$ if and only if there exists exactly one cycle $c_0 \in \cnh{W}$ together with its cyclic permutations such that $f(c_0) \prec f_0(c_0)$ $($under division of $\mathbb{Z}^{+})$, and $f(c) = f_0(c)$ for any other $c \in \cnh{V}$.
\end{lem}

\begin{prf}
It is routine to check the sufficiency part. To prove the necessity part, we suppose that $(H, W, f_0)\prec (H, W, f)$, surely we have $f_0 \neq f$. If $c_1, c_2 \in \cnh{W}$ are not cyclic permutations, $f_0(c_1) \neq f(c_1)$ and $f_0(c_2) \neq f(c_2)$, then defining a cycle function $f'$ with respect to $H$ and $W$ by $f'(c_1) = f(c_1)$ and $f'(c) = f_0(c)$ for any other $c \in \cnh{V}$ (certainly including $c_2$), we have $(H, W, f_0) \lneq (H, W, f') \lneq (H, W, f)$ since $f_0, f'$ are different at $c_1$ and $f', f$ are different at $c_2$. This contradicts $(H, W, f_0) \prec (H, W, f)$. So there exists exactly one cycle $c_0 \in \cnh{W}$ together its cyclic permutations such that $f_0(c_0) \neq f(c_0)$. Certainly, $f(c_0)\in \mathbb{Z}^+$ since otherwise $f_0 (c_0)= f(c_0)=\infty$ so that $f_0 = f$. Furthermore, we also have $f_0 (c_0)\in \mathbb{Z}^+$ since otherwise $f_0 (c_0)= \infty$, and by defining a cycle function $g_1$ with respect to $H$ and $W$ by $g_1(c_0)= 2f(c_0)$ and $g_1(c)=f(c)$ for any other $c \in \cnh{V}$, we obtain that $(H, W, f_0)\lneq (H, W , g_1) \lneq (H, W , f)$, a contradiction. If there exists $m \in \mathbb{Z}^+$ such that $f(c_0) \mid m \mid f_0(c_0)$ and $f_0(c_0) \neq m$, then defining a cycle function with respect to $H$ and $W$ by $g_2(c_0) = m$ and $g_2(c) = f(c)$ for any other $c \in \cnh{V}$, we have $(H, W, f_0) \lneq (H, W, g_2) \leq (H, W, f)$ which leads to $g_2 = f$ so that $m = f(c_0)$. Thus we obtain that $f(c_0) \prec f_0(c_0)$.
\end{prf}

We also need

\begin{lem} \label{Case111}
For any $(H_0, W_0, f_0), (H, W, f)\in \mathcal{CT}(\G)$ with $(H_0, W_0, f_0) \leq (H, W, f)$, if $H_0 \subsetneq H$, $W_0 \setminus H = W$, $f_0(c) = f(c)$ for all $c \in C(W_0)$ and no vertex in $H \setminus (H_0 \cup W_0)$ has index $1$ relative to $H_0$, then for any $(H_0, W', f') \in \mathcal{CT}(\G)$ with $(H_0, W_0, f_0) \leq (H_0, W', f') \lneq (H, W, f)$, one has $W_0 = W'$ and $f_0 = f'$.
\end{lem}

\begin{prf}
It follows from $W_0 \setminus H = W$ that $W_0 = W \cup (W_0 \cap H)$. Moreover, we observe that $W_0 \subseteq W'$ and that $W' \setminus H \subseteq W$ which means that any vertex in $W'$ lies either in $W$ or in $H \setminus H_0$. Because there exists no vertex in $H \setminus (H_0 \cup W_0)$ which has index $1$ relative to $H_0$, we know that each vertex in $W'$ also lies in $W_0$. So we obtain that $W'= W_0$. Then noticing that $f(c) \mid f'(c) \mid f_0(c)$ and $f(c) = f_0(c)$ for any $c \in \cnh{W_0}$, we have $f_0 = f'$.
\end{prf}
\medskip

Now we determine the join and meet of any two elements in $\mathcal{CT}(\G)$. Let $T_1=(H_1, W_1, f_1)$ and $T_2=(H_2, W_2, f_2)$ be arbitrary congruence triples of $\G$. Seeing Figure~2.1, we denote some subsets of $V$. Obviously, both $H_1\cap H_2$ and $H_1 \cup H_2$ are hereditary subsets. Note that vertices in $W_1\cup W_2$ have index at most 1 relative to $H_1 \cup H_2$. Denote the set of vertices in $(W_1 \cup W_2) \setminus (H_1 \cup H_2)$ which have index 0 relative to $H_1 \cup H_2$ by $V_0$, the set $(W_1 \cap W_2) \setminus V_0$ by $X$ (which consists of vertices with index 1 relative to both $H_1 \cup H_2$ and $H_1 \cap H_2$) and the set
$$
\{v\in{(W_1\cup W_2) \setminus (H_1\cup H_2)}: (\exists \alpha \in \pth(\G)) \;v=s(\alpha), u(\alpha) \subseteq W_1 \cup W_2, r(\alpha)\in V_0\}
$$
by $J$. Certainly, $V_0 \subseteq J$ since we consider a vertex as a trivial path.

\begin{figure}[thp]
\centering
	\begin{tikzpicture}[scale=0.5]
	\tikzstyle{every node}=[font=\small,scale=0.9]
	
	\draw  (-1.5,-1.5) circle (1.95);
	\draw (1.5,-1.5) circle (1.95);
	\begin{scope}[fill=red!20,fill  opacity=0.5]
	\draw [clip] (-1.5,-1.5) circle (1.95);
	\fill (1.5,-1.5) circle (1.95);
	\end{scope}
	
	\begin{scope}[fill=blue!20,fill  opacity=0.5]
	\draw [clip] (-1.5,-1.5) circle (1.95);
	\fill (-2,2) circle (3);
	\end{scope}
	
	\begin{scope}[fill=yellow!20,fill  opacity=0.5]
	\draw [clip] (-2,2) circle (3);
	\fill  (2,2) circle (3);
	\end{scope}
	
	\begin{scope}[fill=blue!20,fill  opacity=0.5]
	\draw [clip] (2,2) circle (3);
	\fill (1.5,-1.5) circle (1.95);
	\end{scope}
	
	\begin{scope}[fill=blue!40,fill  opacity=0.5]
	\draw [clip] (2,2) circle (3);
	\draw[clip]  (-2,2) circle (3);
	\fill   (-1,1.5)  rectangle   (1,5) ;
	\end{scope}

	\begin{scope}[fill=brown!,fill  opacity=0.3]
	\draw [clip] (2,2) circle (3) ;
	\fill (-6,1.5)  rectangle   (6,2.8);
	
	\end{scope}
	\begin{scope}[fill=brown!,fill  opacity=0.3]
	\draw [clip] (-2,2) circle (3);
	\fill (-6,1.5)  rectangle   (6,2.8);
	
	\end{scope}
	
	\begin{scope}[fill=blue!10,fill  opacity=0.3]
	\draw [clip] (2,2) circle (3) ;
	\fill (-6,-1)  rectangle   (6,2.8);
	\end{scope}
	
	\begin{scope}[fill=blue!10,fill  opacity=0.3]
	\draw [clip] (-2,2) circle (3) ;
	\fill (-6,-1)  rectangle   (6,2.8);
	\end{scope}
	
	\node[below] at (-5,5) {$W_1$};
	\node[below] at (5,5) {$W_2$};
	\node[below] at (-3,-3.5) {$H_2$};
	\node[below] at (3,-3.5) {$H_1$};
	\node[above]at (-1.6,-0.7) {$W_1\cap H_2$};
	\node[above]at (1.8,-0.7) {$W_2\cap H_1$};
	\node[above]at (0,-2) {$H_2\cap H_1$};
	\node[above] at (0,2.3) {$X$};
	\node[left] at (3,1) {$V_0$};
	\node[left] at (-3,1.5) {$J$};
	\node[below] at (3,-3.5) {$H_1$};
	
	\end{tikzpicture}

{Figure 2.1}
\end{figure}

\begin{lem} \label{J}
If $H'$ is a hereditary subset of $V$ such that $H_1 \cup H_2 \subseteq H'$ and if any vertex in $(W_1 \cup W_2) \setminus H'$ has index $1$ relative to $H'$, then $J \subseteq H'$.
\end{lem}

\begin{prf}
We observe that $V_0\subseteq H'$ since each vertex in $V_0$ has index 0 relative to $H_1 \cup H_2$ which means from $H_1 \cup H_2 \subseteq H'$ that each vertex in $V_0$ has index 0 relative to $H'$. So we need to prove that $J\setminus V_0\subseteq H'$. For any $v \in {J\setminus V_0}$, we see from the definition of $J$ that there exists a path $\alpha = e_1 e_2 ... e_m$ satisfying $u(\alpha) \subseteq W_1 \cup W_2, s(\alpha)= s(e_1)= v$ and $r(\alpha)= r(e_m)\in V_0$. Then it is clear that $s(e_m)$ has index 0 relative to $H'$ and hence belongs to $H'$. In this way, we inductively see that $v \in H'$ as required.
\end{prf}

\begin{lem} \label{meet}
The meet $T_l = (H_l, W_l, f_l)$ of $T_1$ and $T_2$ can be determined as follows: $H_l = H_1 \cap H_2$, $W_l =(W_1 \cap H_2)\cup (W_2\cap H_1) \cup X$ and $f_l(c) = [f_1(c),f_2(c)]$ for $c \in \cnh{V}$.
\end{lem}

\begin{prf}
To show that $T_l\in \mathcal{CT}(\G)$, we first need to prove that for an arbitrary $v \in (W_2 \cap H_1) \cup (W_1 \cap H_2)$, $v$ has index 1 relative to $H_1 \cap H_2$. If $v \in (W_2 \cap H_1)$, then we observe from $v\in H_1$ that $r(e)\in H_1$ for each $e\in s^{-1} (v)$. Moreover, $v\in W_2$ means that there exists a unique edge $e_v \in s^{-1} (v)$ satisfying $r(e_v)\notin H_2$. So we get $r(e_v)\notin H_1 \cap H_2$, and the ranges of the rest edges in $s^{-1} (v)$ are all in $H_1 \cap H_2$. That is to say, $v$ has index 1 relative to $H_1 \cap H_2$. If $v$ in $W_1 \cap H_2$, then we similarly see that $v$ has index 1 relative to $H_1 \cap H_2$. Furthermore, note that $W_2 \cap H_1$ and $W_1 \cap H_2$ are subsets of hereditary subsets which leads to $\cnh{W_l} = \cnh{X} \cup \cnh{W_2 \cap H_1} \cup \cnh{W_1\cap H_2}$ and that no cycles has vertices lying in $(W_1 \cap W_2) \cap V_0$. Obviously, for each $c \in \cnh{H_1 \cap H_2}$, $f_l(c) = 1$. For any $c \in \cnh{V} \setminus (\cnh{W_l} \cup \cnh{H_1 \cap H_2})$, we have either $f_1(c) = \infty$ or $f_2(c) = \infty$ which leads to $f_l(c) = \infty$.

Next, we prove that $T_l$ is a lower bound of $T_1$ and $T_2$. Easy to see that $H_1 \cap H_2 \subseteq H_1$, $W_l\setminus H_1= X \cup(W_1 \cap H_2)\subseteq W_1$. For all $c \in \cnh{V}$, we certainly have $f_1 (c)\mid f_l (c)$. So we obtain that $T_l \leq T_1$. Similarly, $T_l \leq T_2$. Hence $T_l$ is a lower bound of $T_1$ and $T_2$.

Finally, suppose that $(H', W', f')$ is an arbitrary lower bound of $T_1$ and $T_2$. To see that $(H', W', f') \leq T_l$, we only need to show that $W' \setminus (H_1 \cap H_2) \subseteq W_l$ since the other required conditions are obvious. In fact, we observe from $H'\subseteq H_1 \cap H_2$, $W'\setminus H_1 \subseteq W_1$ and $W'\setminus H_2 \subseteq W_2$ that $W'\setminus (H_1 \cap H_2) \subseteq (W_1\cap W_2)\cup (W_2\cap H_1) \cup (W_1 \cap H_2)$. Moreover, for every vertex $v$ in $W'\cap ({W_1\cap W_2})$, we know that $v$ has index 1 or 2 relative to $H_1 \cap H_2$. Suppose that there exists a vertex $v_0 \in W'\cap ({W_1\cap W_2})$ which has index 2 relative $H_1 \cap H_2$. Then $v_0$ must have index at least 2 relative to $H'$ which contradicts the fact that $v_0\in W'$. So we get $W'\cap ({W_1\cap W_2}) \subseteq X$ as required.

In summary, $T_l$ is the meet of $T_1$ and $T_2$.
\end{prf}

\begin{lem} \label{join}
The join $T_u = (H_u, W_u, f_u)$ of $T_1$ and $T_2$ can be determined as follows: $H_u = H_1\cup H_2 \cup J$, $W_u = (W_1 \cup W_2) \backslash{H_u}$ and $f_u(c) = (f_1(c), f_2(c))$ for $c \in \cnh{V}$.
\end{lem}

\begin{prf}
First if $e \in E$ and $s(e) \in J$, then from the definition of $J$, there exists $\al \in \pth(\G)$ such that $u(\al) \subseteq W_1 \cup W_2$, $s(\al) = s(e)$ and $r(\al) \in V_0$. We may suppose that $\al$ is not a vertex since otherwise $s(e)$ has index 0 relative to $H_1 \cup H_2$ which means that $r(e) \in H_1 \cup H_2$. If $e$ is a prefix of $\al$, then we clearly have $r(e) \in J$. If $e$ is not a prefix of $\al$, then noticing that the range of the first edge in $\al$ is not in $H_1 \cup H_2$ and $s(e)$ has index 1 relative to $H_1 \cup H_2$, we must have $r(e) \in H_1 \cup H_2$. So we prove that $H_u$ is a hereditary subset. Now for any $v$ in $(W_1\cup W_2)\setminus H_u$, noticing the definition of $J$ and that vertices in $(W_1 \cup W_2) \setminus (V_0 \cup H_1 \cup H_2)$ have index 1 relative to $H_1 \cup H_2$, $v$ has index 1 relative to $H_u$. By the definition of $J$, we see that there does not exist cycle $c$ such that $u(c) \subseteq \cnh{J}$. So if $c \in \cnh{H_u} = \cnh{H_1 \cup H_2}$, we have $f_u(c) = 1$. Moreover, if $c \in \cnh{V} \setminus \cnh{W_u \cup H_u}$, then $f_u(c) = \infty$ since $f_1(c) = f_2(c) = \infty$. Thus we obtain that $T_u \in \mathcal{CT} (\G)$.

It is clear that $T_u$ is an upper bound of $T_1$ and $T_2$. Now for an arbitrary upper bound $(H', W', f')$ of $T_1$ and $T_2$, we need to check that $H_u \subseteq H'$ and $W_u \setminus H' \subseteq W'$ since we obviously have $f'(c) \mid f_u(c)$ for any $c \in \cnh{V}$. It follows from $W_1\setminus H'\subseteq W'$ and $W_2\setminus H'\subseteq W'$ that $(W_1\cup W_2)\setminus H'\subseteq W'$. That is, any $v \in W_1\cup W_2$ is either in $H'$ or in $W'$. It follows from Lemma~\ref{J} that $H_u \subseteq H'$ so that $W_u \setminus H' = (W_1\cup W_2 )\setminus H'\subseteq W'$.

In conclusion, $T_u$ is the join of $T_1$ and $T_2$.
\end{prf}

\begin{cor}\label{distributive}
Given a fixed hereditary subset $H$ of $V$, the subset $L = \{(H, W, f) : (H, W, f) \in \mathcal{CT}(\G)\}$ is a distributive sublattice of $\mathcal{CT}(\G)$.
\end{cor}

\begin{prf}
For any $(H, W_3, f_3), (H, W_4, f_4) \in L$, we see from Lemmas~\ref{meet} and \ref{join} that $(H, W_3, f_3)\wedge (H, W_4, f_4)= (H, W_3 \cap W_4, f_l)$, where $f_l(c)= [f_3 (c),f_4(c)]$ for $c \in \cnh{V}$ and $(H, W_3, f_3) \vee (H, W_4, f_4)= (H, W_3 \cup W_4, f_u)$, where $f_u(c) = (f_3 (c),f_4(c))$ for $c \in \cnh{V}$. So $(L, \leq)$ is a sublattice of $(\mathcal{CT}(\G), \leq)$. Moreover, given another congruence triple $(H, W_5, f_5)$ in $L$, we have $[(f_3(c),f_4(c)),f_5 (c)] = ([f_3(c),f_5(c)], [f_4(c), f_5(c)])$ for any $c \in \cnh{V}$ (see for instance \cite{Fofanova}) and that $(W_3 \cup W_4) \cap W_5 =(W_3 \cap W_5)\cup(W_4\cap W_5)$. Therefore, $(L, \leq)$ is a distributive lattice.
\end{prf}

\begin{eg}
For the graph given in Figure~2.2, we take hereditary subsets $H_1 = \{v_3, v_4\}$, $H_2 = \{v_9\}$ and subsets $W_1 = \{v_1, v_2, v_5, v_7, v_8\}$, $W_2 = \{v_1, v_4, v_6, v_8, v_{10}, v_{11}\}$ of vertices which have index 1 respectively relative to $H_1$ and $H_2$. Then obviously, $H_1 \cap H_2 = \emptyset = H_2 \cap W_1$, $H_1 \cap W_2 = \{v_4\}$, $W_1 \cap W_2 = \{v_1, v_8\}$ and $H_1 \cup H_2 = \{v_3, v_4, v_9\}$. By the corresponding definitions, we see that $V_0 = \{v_6, v_7, v_8\}$ and $J = \{v_5, v_6, v_7, v_8, v_{10}, v_{11}\}$. Note that $\cnh{V} = \{e_1e_2, e_2e_1, e_3e_4, e_4e_3\}$. According to the definition of an (extended) cycle function, we define $f_1(e_1e_2) = f_1(e_2e_1) = 6$, $f_1(e_3e_4) = f_1(e_4e_3) = 1$ and $f_2(e_1e_2) = f_2(e_2e_1) = f_2(e_3e_4) = f_2(e_4e_3) = \infty$. Then by Lemmas~\ref{meet} and \ref{join}, we have $H_l = \emptyset, W_l = \{v_1, v_4\}, f_l = f_2$ and $H_u = \{v_3, v_4, v_5, v_6, v_7, v_8, v_9, v_{10}, v_{11}\}, W_u = \{v_1, v_2\}, f_u = f_1$.
\end{eg}

\medskip
\begin{center}
\begin{tikzpicture}
\node at (0, 0) (v1) {$\bullet$};
\node at (1.5, 0) (v2) {$\bullet$};
\node at (1.5, -1.5) (v3) {$\bullet$};
\node at (3, 0) (v4) {$\bullet$};
\node at (3, -3) (v5) {$\bullet$};
\node at (4.5, -1.5) (v6) {$\bullet$};
\node at (4.5, 0) (v7) {$\bullet$};
\node at (4.5, 1.5) (v8) {$\bullet$};
\node at (6, 0) (v9) {$\bullet$};
\node at (6, 3) (v10) {$\bullet$};
\node at (6, -3) (v11) {$\bullet$};
\node at (3, 3) (v12) {$\bullet$};

\draw[->] (v1) node[above] {$v_1$} to [bend left = 30] node[above] {$e_1$} (v2);
\draw[->] (v2) node[above] {$v_2$} to [bend left = 30] node[below] {$e_2$} (v1);
\draw[->] (v2) to (v3); \draw[->] (v6) to (v4);
\draw[->] (v3) node[below] {$v_3$} to [bend left = 20] node[above, sloped] {$e_3$} (v4);
\draw[->] (v4) to [bend left = 20] node[below, sloped] {$e_4$} (v3);
\draw[->] (v2) to (v4) node[above] {$v_4$};
\draw[->] (v5) node[below] {$v_5$} to (v4);
\draw[->] (v5) to (v6) node[below] {$v_6$};
\draw[->] (v8) to (v9); \draw[->] (v10) to (v9);
\draw[->] (v7) node[above] {$v_7$} to (v4);
\draw[->] (v7) to (v9) node[right] {$v_9$};
\draw[->] (v6) to [bend left = 20] (v9);
\draw[->] (v6) to [bend right = 20] (v9);
\draw[->] (v8) node[above] {$v_8$} to (v4);
\draw[->] (v10) node[above] {$v_{10}$} to (v8);
\draw[->] (v12) node[above] {$v_{12}$} to (v8);
\draw[->] (v11) node[below] {$v_{11}$} to (v6);
\end{tikzpicture}

Figure~2.2
\end{center}

We continue to prove Theorem~\ref{2 upper semimodular}. Suppose that $T_1\succ T_l$ and that $T_2\succ T_l$. If $H_1= H_2$, then we see from Corollary \ref{distributive} that $T_u \succ T_1, T_2$ since $H_u = H_l = H_1 = H_2$. Thus, we only need to prove the cases where $H_1 \neq H_2$.
\medskip\\
\textbf{Case 1.} $H_1\subsetneq H_2$,

By Lemma~\ref{meet} we see that $H_l = H_1$ and $W_l = (W_1 \cap H_2) \cup X$ (Note that $\cnh{W_l} = \cnh{W_1 \cap H_2} \cup \cnh{X}$ since $H_2$ is hereditary). Then noticing that $(H_l, W_l, f_l) \prec (H_2, W_2, f_2)$ we obtain from Lemma \ref{cover of H}(1) that $W_l \setminus H_2 = X = W_2$ (which means that $X = W_1 \cap W_2 = W_2$), $f_l(c)=f_2(c)$ for any $c \in \cnh{W_l}$. Moreover, combining with $(H_l, W_l, f_l) \prec (H_1, W_1, f_1)$ we have $W_2 \subseteq W_l \subseteq W_1$.\\
\textbf{Subcase 1.1} $W_l \subsetneq W_1$.

We see from the above paragraph and Lemma~\ref{cover of W} that $|W_1 \setminus W_l| = 1$ and $f_l = f_1$ (that is, $f_2(c) \mid f_1 (c)$ for all $c\in \cnh{V}$). Notice that $f_1 (c)= \infty$ for $c\in \cnh{W_1} \setminus \cnh{W_l}$ since $f_2(c)= \infty$ for those $c$. Hence we have $f_1(c) = f_l(c) = f_2(c)$ for any $c \in \cnh{W_1}$. \\
\textbf{Subcase 1.1.1} $W_1 \setminus W_l \nsubseteq V_0$ which means that $V_0 = \emptyset$.

It follows from Lemma~\ref{join} that $H_u = H_2$, $W_u = W_1 \setminus H_2$ and $f_u = f_2$. By Lemma~\ref{cover of W}, we immediately have $T_2 \prec T_u$. Let $(H_1, W_1, f_1)\lneq (H', W', f')\leq (H_u, W_u, f_u)$. We know from Lemmas~\ref{cover of H}(2) and \ref{Case111} that $H_1 \subsetneq H'\subseteq H_2$. It follows from ${(W_1 \cap H_2)\setminus H'} \subseteq  W_1 \setminus H' \subseteq W'$ that every vertex in $(W_1 \cap H_2)\setminus H'$ has index $1$ relative to $H'$. Besides, we observe that each vertex in $W_2$ has index $1$ relative to both $H_1$ and $H_2$ and hence it also has index $1$ relative to $H'$. Defining a cycle function $f''$ with respect to $H'$ and $W_2 \cup ((W_1 \cap H_2) \setminus H')$ by $f''(c) = f_2(c)$ for any $c \in \cnh{W_2 \cup ((W_1 \cap H_2) \setminus H')}$, we see that $(H', W_2 \cup ((W_1 \cap H_2) \setminus H'), f'') \in \mathcal{CT}(\G)$ and that $T_l = (H_1, W_2 \cup (W_1 \cap H_2), f_l) \lneq (H', W_2 \cup ((W_1 \cap H_2) \setminus H'), f'') \leq T_2$ which leads to $H'= H_2$. Thus it follows from $W_1\setminus H_2 \subseteq W'$ and $W'\subseteq W_1 \setminus H_2$ that $W'= W_1 \setminus H_2$. Again noticing that $f_u(c) \mid f'(c) \mid f_1(c), f_1(c) = f_u(c)$ for all $c$ in $\cnh{W_1 \setminus H_2}$, we obtain that $f' = f_u$ so that $(H', W', f')= (H_u, W_u, f_u)$. Hence, we have $T_1 \prec T_u$. \\
\textbf{Subcase 1.1.2} $W_1 \setminus W_l \subseteq V_0$ which means that $W_1 \setminus W_l = V_0$.

Recall that $J = \{v\in W_1 \setminus H_2: (\exists \alpha \in \mbox{\pth}(\G)) ~v= s(\alpha), u(\alpha)\subseteq W_1, r(\alpha) \in V_0\}$ in this case. We easily observe that all vertices passing by any cycle in $W_1$ must lie entirely in $W_2 \setminus J$ or entirely in $W_1 \cap H_2$, which means that $\cnh{W_1}= \cnh{W_2 \setminus J} \cup \cnh{W_1 \cap H_2} = \cnh{W_2} \cup \cnh{W_1 \cap H_2}$. Then it follows from Lemma~\ref{join} that $H_u = H_2 \cup J$, $W_u = W_2 \setminus J$ and $f_u = f_2$.

Let $(H_1, W_1, f_1)\lneq (H', W', f')\leq (H_2 \cup J, W_2 \setminus J, f_u)$. We obtain from Lemmas~\ref{cover of H}(2) and \ref{Case111} that $H_1 \subsetneq H'$. Note that $H' \subseteq H_2 \cup J$. If $J \cap H' \neq \emptyset$, then we know from the definition of $J$ that $V_0 \subseteq H'$. Noticing the unique element in $V_0$ has index 0 relative to $H_2$ and has index 1 relative to $H_1$, in this case we have $H' \cap (H_2 \setminus H_1) \neq \emptyset$. If $J \cap H' = \emptyset$, then we also have $H' \cap (H_2 \setminus H_1) \neq \emptyset$ since $H_1 \subsetneq H'$. So we get $H_1 \subsetneq H' \cap H_2$. Similar to the discussion in Case~1.1.1, if we define a cycle function $f''$ with respect to $H' \cap H_2$ and $W_2 \cup ((W_1 \cap H_2) \setminus H')$ by $f''(c) = f_2(c)$ for any $c \in \cnh{W_2 \cup ((W_1 \cap H_2) \setminus H')}$, then we obtain that $(H' \cap H_2, W_2 \cup ((W_1 \cap H_2) \setminus H'), f'') \in \mathcal{CT}(\G)$ and that $T_l \lneq (H' \cap H_2, W_2 \cup ((W_1 \cap H_2) \setminus H'), f'') \leq T_2$, which means that $H' \cap H_2 = H_2$. So we get $H_2 \subseteq H'$. Thus, we know from $(H_1, W_1, f_1) \leq (H', W', f')$ that every vertex in $W_1$ is either in $H'$, or in $W'$. It follows from Lemma~\ref{J} that $J \subseteq H'$ so that $H'= H_2 \cup J$. From $W_2 \setminus J = W_1 \setminus H' \subseteq W' \subseteq W_2 \setminus J$ we know that $W'= W_2 \setminus J$. Moreover, since $f_1(c) = f_2(c)= f_u(c)$ for all $c$ in $\cnh {W_1}$, and $f_u(c) \mid f'(c) \mid f_1 (c)$ for all $c \in \cnh{W_2\setminus J}$, we have $f'(c) = f_u(c)$ for all $c \in \cnh{W_2\setminus J}$. So $f'= f_u$ for all $c\in \cnh{V}$, that is, $(H', W', f')= (H_2 \cup J, W_2 \setminus J, f_u)$ which leads to $T_1 \prec T_u$.

Let $(H_2, W_2, f_2)\lneq (H'', W'', f'')\leq (H_2 \cup J,  W_2 \setminus J, f_u)$. It follows from Lemma~\ref{Case111} that $H_2 \subsetneq H''$. Then we see from $H''\cap J \neq\emptyset$ that $V_0 \subseteq H''$, certainly, $J\subseteq H''$. So we have $H''= H_2 \cup J, W'' = W_2 \setminus J$. Noticing that $f_2(c) = f_u(c)$, we observe that $(H'', W'', f'')=(H_2 \cup J, W_2 \setminus J, f_u)$ which also leads to $T_2 \prec T_u$.\\
\textbf{Subcase 1.2} $W_l = W_1$.

Recall from the hypotheses in this case, Lemmas~\ref{meet}, \ref{join} and \ref{cover of H}(1) that $T_l = (H_1, W_1, f_l)$, where $f_l(c) = f_2(c)$ for any $c \in \cnh{W_l}=\cnh{W_1}$, and $T_u = (H_2, W_2, f_u)$ since $V_0 = \emptyset$, where $f_u(c) = f_1(c)$ for any $c \in \cnh{W_1}$. Notice that $f_1(c)= 1$ for any $c\in \cnh{W_1 \cap H_2}$ since $f_1(c)\mid f_2(c)$ and $f_l(c) = f_2(c)=1$. We obtain that $f_1(c)= f_l(c)= 1$ for all $c\in \cnh{W_1 \cap H_2}$. Thus we can see from Lemma~\ref{samehw} that there exists exactly one cycle $c_0 \in \cnh{W_2}$ together with its cyclic permutations such that $f_u (c_0)= f_1(c_0) \prec f_l(c_0) = f_2(c_0)$ and $f_1(c) = f_l(c)$ for any other $c \in \cnh{V}$. It follows again from Lemma~\ref{samehw} that $T_2 \prec T_u$.

If $(H_1, W_1, f_1) \lneq (H', W', f') \leq (H_2, W_2, f_u)$, then we see from Lemmas~\ref{cover of H}(2) and \ref{Case111} that $H_1 \subsetneq H'$. Defining a cycle function $f''$ with respect to $H'$ and $W_2 \cup ((W_1 \cap H_2) \setminus H')$ by $f''(c) = f_2(c)$ for any $c \in \cnh{W_2 \cup ((W_1 \cap H_2) \setminus H')}$, we see that $(H', W_2 \cup ((W_1 \cap H_2) \setminus H'), f'') \in \mathcal{CT}(\G)$ and that
$T_l = (H_1, W_1, f_l) \lneq (H', W_2 \cup ((W_1 \cap H_2) \setminus H'), f'') \leq T_2$,
which leads to $H'= H_2$. Thus it follows from $W_1\setminus H_2 \subseteq W'$ and $W'\subseteq W_1 \setminus H_2$ that $W'= W_1 \setminus H_2 = W_2$. Again noticing that $f_u(c) \mid f'(c) \mid f_1(c), f_1(c) = f_u(c)$ for all $c$ in $\cnh{W_2}$, so we obtain that $f' = f_u$ and $(H', W', f')= (H_u, W_u, f_u)$. Hence, we have $T_1 \prec T_u$.  \medskip \\
\textbf{Case 2.} $H_2\subsetneq H_1$.

Similar to case 1. \medskip \\
\textbf{Case 3.} $H_1\nsubseteq H_2, H_2 \nsubseteq H_1$.

We see from Lemma~\ref{meet} that $H_l \subsetneq H_1, H_l \subsetneq H_2$ and that $W_l = (W_2 \cap H_1)\cup (W_1 \cap H_2) \cup X$. It follows from Lemma~\ref{cover of H}(1) that $W_l \setminus H_1 = W_1, W_l \setminus H_2 = W_2$, $f_l(c)=f_1(c)= f_2 (c)$ for any $c\in \cnh{W_l}$ (which means by (\ref{extending}) that $f_1(c) = 1$ for $c \in \cnh{W_1 \cap H_2}$ and that $f_2(c) = 1$ for $c \in \cnh{W_2 \cap H_1}$). We observe that $X = W_1 \cap W_2$, $W_l = W_1 \cup W_2$ and $V_0 = \emptyset$. Then we obtain from Lemma~\ref{join} that $(H_u, W_u, f_u)= (H_1\cup H_2, W_1 \cap W_2, f_u)$, where $f_u (c)= f_1 (c)= f_2 (c)$ for all $c\in \cnh{W_u}$.

Let $(H_1, W_1, f_1)\lneq (H', W', f')\leq (H_1\cup H_2, W_1 \cap W_2, f_u)$. Again we see from Lemmas~\ref{cover of H}(2) and \ref{Case111} that $H_1 \subsetneq H' \subseteq H_1 \cup H_2$ which leads to $H_1 \cap H_2 \subsetneq H' \cap H_2$. Noticing that $W_1 \setminus H' \subseteq W'$ and that every vertex in $(W_1 \cap H_2)\setminus H'$ has index $1$ relative to $H' \cap H_2$, we observe that $(H' \cap H_2, W_l\setminus (H' \cap H_2), f'') \in \mathcal{CT}(\G)$, where $f''(c) = f_l(c)$ for any $c \in \cnh{W_l\setminus (H' \cap H_2)}$.
Moreover, $T_l \lneq (H' \cap H_2, W_l\setminus (H' \cap H_2), f'') \leq T_2$. Hence, we get $H' \cap H_2 = H_2$. Combining this with $H_1 \subsetneq H' \subseteq H_1 \cup H_2$, we get $H'= H_1\cup H_2$. Furthermore, $W_1\setminus (H_1 \cup H_2) = (W_1 \cap W_2)\subseteq W' \subseteq (W_1 \cap W_2)$ which means that $W'= W_1 \cap W_2$. Now it is easy to see that $(H', W', f')= T_u$ and hence $T_1 \prec T_u$.
Similarly, $T_2 \prec T_u$.\medskip

As a consequence, we obtain that $(\mathcal{CT}(\G), \leq)$ is upper semimodular. By Proposition~\ref{isom}, we know that $(\mathcal{CS}, \subseteq)$ is upper semimodular. However, they are not lower semimodular as the following example illustrates.

\begin{eg}
Let $\G$ be the graph where $V= \{v_0, v_1, v_2\}$, $E= \{e_1, e_2\}$, $s(e_1) = s(e_2) = v_0$, $r(e_1) = v_1$ and $r(e_2) = v_2$.
We see that $H_1=\{v_1\}, H_2=\{v_2\}$ are two hereditary subsets. $v_0$ has index $1$ relative to both $H_1$ and $H_2$. Let $W_1=W_2 =\{v_0\}$. Then we clearly have $(H_1, W_1, \emptyset), (H_2, W_2, \emptyset) \in \mathcal{CT}(\G)$ and $V_0 = \{v_0\}$. It follows from Lemma \ref{join} that $(H_1, W_1, f_1)\vee (H_2, W_2, f_2)=(V, \emptyset, \emptyset)$.

We claim that $(H_1, W_1, \emptyset)$ and $(H_2, W_2, \emptyset)$ are both covered by $(V, \emptyset, \emptyset)$. In fact, suppose $(H', W', \emptyset)\in \mathcal{CT}(\G)$ satisfying $(H_1, W_1, \emptyset) \lneq (H', W', \emptyset) \leq (V, \emptyset, \emptyset)$. If $H'= H_1$, then we clearly have $W'=W_1=\{v_0\}$ which leads to $(H', W', f')=(H_1, W_1, f_1)$, a contradiction. If $H'= H_1\cup H_2$, then $W'= \emptyset$ which contradicts the condition that $W_1\setminus (H_1\cup H_2)=\{v_0\}\subseteq W'$. So we get $H'=V$ and $W'=\emptyset$ which lead to $(H_1, W_1, \emptyset) \prec(V, \emptyset ,\emptyset)$. Similarly $(H_2, W_2, \emptyset) \prec (V, \emptyset, \emptyset)$.

By Lemma \ref{meet}, we know that $(H_1, W_1, f_1)\wedge (H_2, W_2, f_2) = (\emptyset, \emptyset, \emptyset)$.  But there exists $(H_1, \emptyset, \emptyset)\in \mathcal{CT}(\G)$ satisfying
$(\emptyset, \emptyset, \emptyset)\lneq (H_1, \emptyset, \emptyset)\lneq (H_1, W_1, \emptyset)$, which means that $(\mathcal{CT}(\G), \leq )$ is not lower semimodular.

In fact, we can also see that the lattice of congruences on the graph inverse semigroup of this graph has the Hasse diagram as in Figure~2.3, where $T = (\{v_1, v_2\}, \emptyset, \emptyset)$.
\end{eg}

\begin{center}\setlength{\unitlength}{1mm}
\begin{picture}(20,46)
\put(10,10){\circle*{1.5}} \put(10,30){\circle*{1.5}}
\put(0,20){\circle*{1.5}} \put(20,20){\circle*{1.5}}
\put(0,30){\circle*{1.5}} \put(20,30){\circle*{1.5}}
\put(10,40){\circle*{1.5}}
\put(10,10){\line(-1,1){10}} \put(10,10){\line(1,1){10}}
\put(0,20){\line(0,1){10}} \put(0,20){\line(1,1){10}}
\put(20,20){\line(-1,1){10}} \put(20,20){\line(0,1){10}}
\put(10,40){\line(-1,-1){10}} \put(10,40){\line(1,-1){10}}
\put(10,40){\line(0,-1){10}}
\put(5,5){$(\emptyset, \emptyset, \emptyset)$}
\put(21,19){$(\{v_2\}, \emptyset, \emptyset)$}
\put(21,29){$(\{v_2\}, \{v_0\}, \emptyset)$}
\put(-19,19){$(\{v_1\}, \emptyset, \emptyset)$}
\put(-24.5,29){$(\{v_1\}, \{v_0\}, \emptyset)$}
\put(5,42){$(V, \emptyset, \emptyset)$}
\put(8.5,25){$T$}
\put(0, 0){Figure 2.3}
\end{picture}
\end{center}

\medskip

\noindent {\bf Acknowledgements} The second author thanks John Meakin for discussions concerning various aspects of this paper.


\end{document}